# Sample autocovariances of long-memory time series

LAJOS HORVÁTH[1] and PIOTR KOKOSZKA[2]

[1]*Department of Mathematics, University of Utah, 155 South 1440 East, Salt Lake City, UT 84112-0090, USA. E-mail: horvath@math.utah.edu*
[2]*Department of Mathematics and Statistics, Utah State University, 3900 Old Main Hill, Logan, UT 84322-3900, USA. E-mail: Piotr.Kokoszka@usu.edu*

We find the asymptotic distribution of the sample autocovariances of long-memory processes in cases of finite and infinite fourth moment. Depending on the interplay of assumptions on moments and the intensity of dependence, there are three types of convergence rates and limit distributions. In particular, a normal approximation with the standard rate does not always hold in practically relevant cases.

## 1. Introduction

Over the last twenty years, long-memory time series have become an important modeling tool in geophysical sciences and also in engineering, computer networks and econometrics.

Asymptotics for the sample autocovariances of long-range dependent linear processes with tail index $\alpha$ are known only in the case $1 < \alpha < 2$; see Kokoszka and Taqqu (1996). In the case $2 < \alpha < 4$, these asymptotics are known only for linear processes with absolutely summable coefficients $\psi(j)$; see Davis and Resnick (1986). To illustrate what kind of results we are interested in, consider the sample variance $\hat{\gamma}_0 = N^{-1} \sum_{t=1}^N X_t^2$. Theorem 2.2 of Davis and Resnick (1986) implies that

$$N a_N^{-2}(\hat{\gamma}_0 - \gamma_0) \overset{d}{\to} \sum_{j=0}^\infty \psi^2(j)\left(S - \frac{\alpha}{\alpha-2}\right), \tag{1.1}$$

where $a_N$ is roughly of the order $N^{1/\alpha}$ and $S$ is an $(\alpha/2)$-stable random variable. Note that the right-hand side of (1.1) involves only the sum of the squared coefficients $\psi(j)$, yet this result is known to hold only if $\sum_{j=0}^\infty |\psi(j)| < \infty$. The question is whether (1.1) holds if one assumes only $\sum_{j=0}^\infty \psi^2(j) < \infty$ (and zero mean) and, if not, what the limit of $N a_N^{-2}(\hat{\gamma}_0 - \gamma_0)$ is in this more general case.

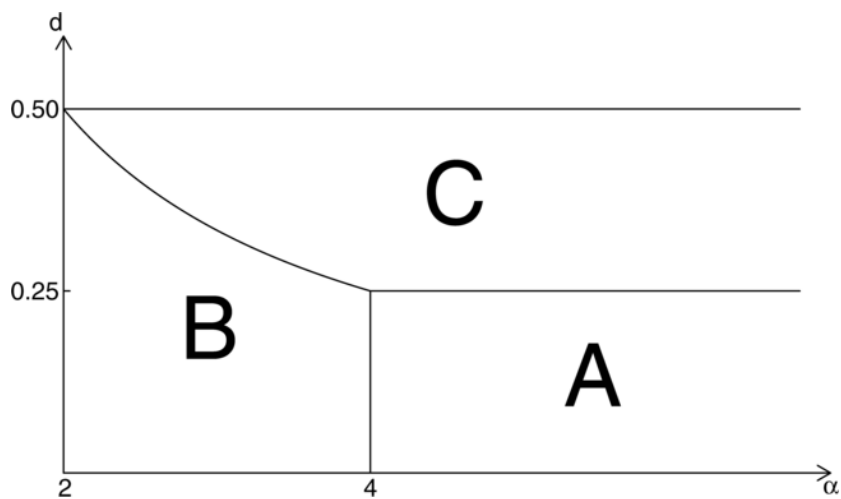

**Figure 1.** Regions with different rates of convergence of sample autocovariances.

We assume that the $\psi(j)$ behave roughly like $j^{d-1}$ with $0 < d < 1/2$. Figure 1 summarizes the convergence rates derived in our paper. The $\alpha$-axis shows the tail index of the innovations, but $\alpha > 4$ is used merely for illustration; our results require only finite fourth moment in that case. The currently known asymptotics correspond to the cases with $d = 0$, except the case when all moments are finite (considered by Hosking (1996)). Our results for the region $4 < \alpha < \infty$ extend the corresponding results of Hosking (1996). In region A, $\hat{\gamma}_0 - \gamma_0$ is of the standard order $N^{-1/2}$; in region B, it is roughly of the order $N^{2/\alpha - 1}$. These two regions show how far the rates established for weakly dependent processes extend. In region C, the asymptotic distribution of $\hat{\gamma}_h - \gamma_h$ does not depend on $h$ and the limit is related to the Rosenblatt process.

The paper is organized as follows. Section 2 introduces the assumptions and the requisite notation. Main results are stated in Section 3. Section 4 contains an important truncation lemma which allows us to assume only the second moment and two lemmas which can be deduced from the work of Surgailis (1982). The proofs of the main results are developed in Section 5.

## 2. Assumptions and notation

We consider the linear process

$$X_t = \sum_{j=0}^{\infty} \psi(j) Z_{t-j}. \tag{2.1}$$

Throughout the paper, we denote by $Z$ a random variable with the distribution of the $Z_t$.

We assume that the $Z_t$ are mean zero i.i.d. and satisfy either $\mathbf{E}Z^4 < \infty$ or for some $2 < \alpha < 4$, $0 \leq p \leq 1$ and a slowly varying function $L(x)$,

$$P[|Z_t| > x] = x^{-\alpha}L(x), \qquad P[Z_t > x]/P[|Z_t| > x] \to p \qquad (x \to \infty). \tag{2.2}$$

In both cases, the second moment is finite and we write

$$\sigma^2 = \mathrm{Var}[Z] = \mathbf{E}Z^2.$$

We assume

$$\sum_{j=0}^{\infty} \psi^2(j) < \infty. \tag{2.3}$$

Assumption (2.3) is a sufficient condition for the $L^2$ and a.s. convergence of the right-hand side of (2.1) if $\mathbf{E}Z_t^2 < \infty$. (Note that we assume $\mathbf{E}Z = 0$.) Condition (2.3) is weaker than the condition $\sum_{j=0}^{\infty} |\psi(j)| < \infty$ and allows us to consider linear processes with long memory.

We focus on the case $\psi(j) \sim C_d j^{d-1}$ and assume

$$\psi(j) = j^{d-1}l(j), \qquad 0 < d < 1/2, \tag{2.4}$$

where $l(\cdot)$ is a function defined for positive real numbers such that $l(u) \to C_d$ as $u \to \infty$. All long memory models in current use, including the fractional ARIMA and the fractional exponential model, satisfy (2.4).

Suppose we observe a realization $X_1, X_2, \ldots, X_{N+H}, N > 1, H \geq 0$. The sample autocovariances are defined as

$$\hat{\gamma}_h = \frac{1}{N}\sum_{t=1}^{N} X_t X_{t+h}, \qquad h = 0, 1, \ldots, H,$$

and the (population) autocovariances are

$$\gamma_h = \mathbf{E}[X_0 X_h] = \sigma^2 \sum_{j=0}^{\infty} \psi(j)\psi(j+h).$$

When working with the innovations satisfying (2.2), we use the norming constants $a_N$ satisfying

$$\lim_{N \to \infty} NP[|Z| > a_N x] = x^{-\alpha}, \qquad x > 0. \tag{2.5}$$

It is well known (CLT for i.i.d. r.v.'s in the stable domain of attraction) that there is an $(\alpha/2)$-stable random variable $S$ such that

$$a_N^{-2}\sum_{t=1}^{N}(Z_t^2 - b_N) \xrightarrow{d} S, \tag{2.6}$$

where $b_N = \mathbf{E}[Z^2 I\{|Z| \le a_N\}]$.

Finally, recall the definition of the Rosenblatt process:

$$U_d(t) = 2\int_{x_1 < x_2 < t}\left[\int_0^t (v-x_1)_+^{d-1}(v-x_2)_+^{d-1}\,\mathrm{d}v\right]W(\mathrm{d}x_1)W(\mathrm{d}x_2),\qquad(2.7)$$

where $W(\cdot)$ is the standard Wiener process on the real line.

## 3. Main results

We want to find the asymptotic distribution of the vector

$$[\hat{\gamma}_h - \gamma_h, h = 0, 1, \ldots, H].$$

We first focus on the case of heavy-tailed innovations satisfying (2.2).

**Theorem 3.1.** *Suppose (2.2) and (2.4) hold.*

(a) *If* $0 < d < 1/\alpha, 2 < \alpha < 4$, *then*

$$Na_N^{-2}[\hat{\gamma}_h - \gamma_h, h = 0, 1, \ldots, H]$$
$$\xrightarrow{d} \left(S - \frac{\alpha}{\alpha - 2}\right)\left[\sum_{j=0}^{\infty}\psi(j)\psi(j+h), h = 0, 1, \ldots, H\right].$$

*(For the above to hold for* $d = 1/4$, *we must additionally assume that* $a_N^{-4}N\ln N \to 0$.)

(b) *If* $1/\alpha < d < 1/2, 2 < \alpha < 4$, *then*

$$N^{1-2d}[\hat{\gamma}_h - \gamma_h, h = 0, 1, \ldots, H] \xrightarrow{d} \sigma^2 C_d^2[U_d(1), h = 0, 1, \ldots, H].$$

**Proof.** A detailed proof with the requisite notation is presented in Section 5. Part (a) follows from Lemmas 5.1 and 5.2, while (b) follows from Lemmas 5.3 and 5.5. □

If $1 < \alpha < 2$, then Theorem 2.1 of Kokoszka and Taqqu (1996) implies that $Na_N^{-2}\hat{\gamma}_h$ converges in distribution, so the autocovariances themselves are not bounded in probability.

We now turn to the case $\mathbf{E}Z^4 < \infty$. If the $\psi(j)$ are absolutely summable, then we have the following well-known result; see Proposition 7.3.3 of Brockwell and Davis (1991).

**Theorem 3.2.** *Suppose* $\sum_{j=0}^{\infty}|\psi(j)| < \infty$ *and* $\mathbf{E}Z^4 = \eta\sigma^4 < \infty$. *Then*

$$N^{1/2}[\hat{\gamma}_h - \gamma_h, h = 0, 1, \ldots, H] \xrightarrow{d} [G_h, h = 0, 1, \ldots, H],$$

*where $[G_h, h = 0, 1, \ldots, H]$ is a mean zero Gaussian vector with*

$$\mathbf{E}[G_h G_{h'}] = (\eta - 3)\gamma_h \gamma_{h'} + \sum_{k=-\infty}^{\infty} [\gamma_k \gamma_{k-h+h'} + \gamma_{k+h'}\gamma_{k-h}]. \tag{3.1}$$

Condition (2.4) implies that the $\psi(j)$ are not absolutely summable. In this case, the following theorem holds.

**Theorem 3.3.** *Suppose $\mathbf{E}Z^4 = \eta\sigma^4 < \infty$ and (2.4) holds.*

(a) *If $0 < d < 1/4$, then*

$$N^{1/2}[\hat{\gamma}_h - \gamma_h, h = 0, 1, \ldots, H] \overset{d}{\to} [G_h, h = 0, 1, \ldots, H],$$

*where $[G_h, h = 0, 1, \ldots, H]$ is a mean zero Gaussian vector with covariances (3.1).*

(b) *If $1/4 < d < 1/2$, then*

$$N^{1-2d}[\hat{\gamma}_h - \gamma_h, h = 0, 1, \ldots, H] \overset{d}{\to} \sigma^2 C_d^2[U_d(1), h = 0, 1, \ldots, H].$$

**Proof.** Use Proposition 5.1 for (a) and Lemmas 5.3 and 5.5 for (b). $\qquad\square$

# 4. Auxiliary lemmas

We begin with a truncation lemma which allows us to extend the results established under the assumption of all finite moments to our setting in which only the second moment is assumed finite.

**Lemma 4.1.** *Suppose $Z$ is a random variable with $\mathbf{E}Z = 0$ and $\mathbf{E}Z^2 < \infty$. There are then bounded random variables $Z(T), T > 0$, such that $\mathbf{E}Z(T) = 0$ and, as $T \to \infty$, $\mathbf{E}[Z(T) - Z]^2 \to 0$ and $\mathbf{E}Z^2(T) \to \mathbf{E}Z^2$.*

**Proof.** Let $\mu(T) = \mathbf{E}[Z\mathbf{I}\{|Z| \le T\}]$. If $\mu(T) = 0$, set $Z(T) = Z\mathbf{I}\{|Z| \le T\}$. If $\mu(T) \ne 0$, then let $\varepsilon(T)$ be a random variable uniform on $[0, |\mu(T)|]$ and independent of $Z$, and set

$$Z(T) = Z\mathbf{I}\{|Z| \le T\} - 2\varepsilon(T)\operatorname{sign}(\mu(T)).$$

Then $\mathbf{E}Z(T) = 0$, by the definition of $Z(T)$. Moreover,

$$\mathbf{E}[Z(T) - Z]^2 = \mathbf{E}[Z\mathbf{I}\{|Z| > T\} - 2\varepsilon(T)\operatorname{sign}(\mu(T))]^2$$
$$\le 2[\mathbf{E}(Z^2\mathbf{I}\{|Z| > T\}) + 4\mathbf{E}\varepsilon^2(T)].$$

Since $\mathbf{E}Z^2 < \infty$, $\mathbf{E}(Z^2\mathbf{I}\{|Z| > T\}) \to 0$, and $\mathbf{E}\varepsilon^2(T) = \mu^3(T)/3 \to 0$ because $\mu(T) \to \mathbf{E}Z = 0$. The assertion $\mathbf{E}Z^2(T) \to \mathbf{E}Z^2$ now also follows. $\qquad\square$

The following lemma is a special case of Lemma 7 of Surgailis (1982).

Denote by $\{W(t)\}$ the standard Wiener process on the real line. Consider the sequence $\{W_N, N \geq 0\}$ of random sequences

$$W_N = \{W_N(k), k = \ldots, -1, 0, 1, \ldots\}.$$

It is convenient to think of $W_N(k)$ as approximately the increment of $W$ over the interval $[k/N, (k+1)/N)$.

**Lemma 4.2.** *Suppose each $W_N(k)$ has all finite moments, $\mathbf{E}W_N(k) = 0$ and $\mathrm{Var}[W_N(k)] = N^{-1}$. Assume that for any real numbers $a < b$,*

$$\sum_{aN \leq k < bN} W_N(k) \xrightarrow{d} W(b) - W(a).$$

*Let $f(x_1, x_2)$ be a measurable function such that*

$$\int_{-\infty}^{\infty} \int_{-\infty}^{\infty} |f(x_1, x_2)|^2 \, \mathrm{d}x_1 \, \mathrm{d}x_2 < \infty.$$

*Suppose $\{f_N, N \geq 1\}$ is a sequence of measurable functions of the form*

$$f_N(x_1, x_2) = \sum_{k \neq k'} c_N(k, k') \mathbf{I}_{[k/N,(k+1)/N)}(x_1) \mathbf{I}_{[k'/N,(k'+1)/N)}(x_2) \tag{4.1}$$

*such that*

$$\int_{-\infty}^{\infty} \int_{-\infty}^{\infty} |f(x_1, x_2) - f_N(x_1, x_2)|^2 \, \mathrm{d}x_1 \, \mathrm{d}x_2 \to 0 \qquad \text{as } N \to \infty. \tag{4.2}$$

*Then, as $N \to \infty$,*

$$\sum_{k \neq k'} c_N(k, k') W_N(k) W_N(k') \xrightarrow{d} \int_{-\infty}^{\infty} \int_{-\infty}^{\infty} f(x_1, x_2) W(\mathrm{d}x_1) W(\mathrm{d}x_2).$$

The next result is similar to Lemma 4 of Surgailis (1982). The extension involves the presence of the lag $h$. To state it, we introduce the coefficients

$$C_{N,h}(k, k') = N^{1-2d} \sum_{t=1}^{N} \psi(t-k) \psi(t+h-k') \tag{4.3}$$

and set

$$f_{N,h}(x_1, x_2) = \sum_{k \neq k' \leq N} C_{N,h}(k, k') \mathbf{I}_{[k/N,(k+1)/N)}(x_1) \mathbf{I}_{[k'/N,(k'+1)/N)}(x_2). \tag{4.4}$$

**Lemma 4.3.** *Consider the function $f$ defined by (5.12) and the sequence of functions $f_N$ defined by (4.4) and (4.3). Suppose (2.4) holds. If $d > 1/4$, then*

$$\int_{-\infty}^{\infty} \int_{-\infty}^{\infty} [f_{N,h}(x_1, x_2) - f(x_1, x_2)]^2 \, \mathrm{d}x_1 \, \mathrm{d}x_2 \to 0 \qquad \text{as } N \to \infty. \tag{4.5}$$

**Proof.** By Lemma 4 of Surgailis (1982), (4.5) holds with $h = 0$. Thus, it is enough to show that

$$\int_{-\infty}^{\infty} \int_{-\infty}^{\infty} [f_{N,h}(x_1, x_2) - f_{N,0}(x_1, x_2)]^2 \, \mathrm{d}x_1 \, \mathrm{d}x_2 \to 0 \qquad \text{as } N \to \infty. \tag{4.6}$$

Relation (4.6) can be verified using the dominated convergence theorem. $\qquad \square$

# 5. Proofs

We begin with the proof of part (a) of Theorem 3.3 as it uses an approach which differs from that of the remaining proofs. The idea is the same as in the proof of Proposition 7.3.3 of Brockwell and Davis (1991), but the summability arguments must be handled with care.

**Proposition 5.1.** *The conclusions of Theorem 3.2 remain true if the assumption $\sum_{j=0}^{\infty} |\psi(j)| < \infty$ is replaced by (2.4) with the restriction $0 < d < 1/4$.*

**Proof.** For sufficiently large $m$, define

$$X_t^{(m)} = \sum_{j=0}^{m} \psi(j) Z_{t-j}, \qquad \hat{\gamma}_h^{(m)} = N^{-1} \sum_{t=1}^{N} X_t^{(m)} X_{t+h}^{(m)}$$

and

$$\gamma_h^{(m)} = \mathbf{E}[X_t^{(m)} X_{t+h}^{(m)}] = \sigma^2 \sum_{j=0}^{m-h} \psi(j) \psi(j+h).$$

Proposition 7.3.2 of Brockwell and Davis (1991) states that, as $N \to \infty$,

$$N^{1/2}[\hat{\gamma}_h^{(m)} - \gamma_h^{(m)}, h = 0, 1, \ldots, H] \xrightarrow{d} [G_h^{(m)}, h = 0, 1, \ldots, H], \tag{5.1}$$

where $[G_h^{(m)}, h = 0, 1, \ldots, H]$ is a mean zero Gaussian vector with

$$\mathbf{E}[G_h^{(m)} G_{h'}^{(m)}] = (\eta - 3) \gamma_h^{(m)} \gamma_{h'}^{(m)} + \sum_{k=-\infty}^{\infty} [\gamma_k^{(m)} \gamma_{k-h+h'}^{(m)} + \gamma_{k+h'}^{(m)} \gamma_{k-h}^{(m)}].$$

The proof is completed by appealing to the standard argument; see Theorem 3.2 on page 28 of Billingsley (1999), which, in addition to (5.1), requires

$$[G_h^{(m)}, h = 0, 1, \ldots, H] \xrightarrow{d} [G_h, h = 0, 1, \ldots, H], \qquad \text{as } m \to \infty \tag{5.2}$$

and, for each $h$ and $\epsilon > 0$,

$$\lim_{m \to \infty} \limsup_{N \to \infty} P\{|N^{1/2}(\hat{\gamma}_h^{(m)} - \gamma_h^{(m)}) - N^{1/2}(\hat{\gamma}_h - \gamma_h)| > \epsilon\} = 0. \tag{5.3}$$

Relation (5.2) follows from the convergence $\gamma_h^{(m)} \to \gamma_h$ for which only the square summability of the $\psi(j)$ is needed.

To prove (5.3), it suffices to show that $\lim_{m \to \infty} \lim_{N \to \infty} N \operatorname{Var}[\hat{\gamma}_h^{(m)} - \hat{\gamma}_h] = 0$, which reduces to verifying that

$$\lim_{N \to \infty} N \operatorname{Var}[\hat{\gamma}_h] = \mathbf{E} G_h^2, \tag{5.4}$$

$$\lim_{m \to \infty} \lim_{N \to \infty} N \operatorname{Var}[\hat{\gamma}_h^{(m)}] = \mathbf{E} G_h^2, \tag{5.5}$$

$$\lim_{m \to \infty} \lim_{N \to \infty} N \operatorname{Cov}(\hat{\gamma}_h^{(m)}, \hat{\gamma}_h) = \mathbf{E} G_h^2. \tag{5.6}$$

To establish (5.4)–(5.6), we need the assumption $d < 1/4$. As in the proof of Proposition 7.3.1 of Brockwell and Davis (1991), we have

$$N \operatorname{Var}[\hat{\gamma}_h]$$
$$= \sum_{|k| < N} \left(1 - \frac{|k|}{N}\right) \left[(\eta - 3)\sigma^4 \sum_i \psi(i)\psi(i+h)\psi(i+k)\psi(i+k+h) + \gamma_k \gamma_k + \gamma_{k+h} \gamma_{k-h}\right].$$

Relation (5.4) now follows from the dominated convergence theorem, for which we need

$$\sum_k \sum_i |\psi(i)\psi(i+h)\psi(i+k)\psi(i+k+h)| < \infty \tag{5.7}$$

and

$$\sum_k \gamma_k^2 < \infty. \tag{5.8}$$

While the square summability of the $\psi(j)$ is sufficient for (5.7) to hold, for (5.8), we need $\sum_{k>0} k^{4d-2} < \infty$, which requires that $d < 1/4$.

The same argument and the convergence $\lim_{m \to \infty} \mathbf{E} G_h^{(m)2} = \mathbf{E} G_h^2$ lead to (5.5). Relation (5.6) follows in a similar manner. □

In the following, we work separately with diagonal and off-diagonal terms.

Fix real numbers $u_h, h = 0, 1, \ldots, H$, and consider the decomposition

$$\sum_{h=0}^{H} u_h(\hat{\gamma}_h - \gamma_h) = D_N + R_N,$$

with the diagonal terms

$$D_N = \sum_{h=0}^{H} u_h d_{N,h}, \qquad d_{N,h} = \frac{1}{N} \sum_{t=1}^{N} \sum_{j=0}^{\infty} \psi(j)\psi(j+h)[Z_{t-j}^2 - \sigma^2],$$

and the off-diagonal terms

$$R_N = \sum_{h=0}^{H} u_h r_{N,h}, \qquad r_{N,h} = \frac{1}{N} \sum_{t=1}^{N} \sum_{i \neq j+h} \psi(j)\psi(i)Z_{t-j}Z_{t+h-i}.$$

Also, define

$$c_j = \sum_{h=0}^{H} u_h c_j(h), \qquad c_j(h) = \psi(j)\psi(j+h). \qquad (5.9)$$

Note that $c_j(h) \sim C_d^2 j^{2d-2}$, so the $c_j$ are absolutely summable.

The next two lemmas are specific to the case of infinite fourth moment and so are established first.

**Lemma 5.1.** *Suppose* $2 < \alpha < 4$ *and that (2.2) and (2.4) hold. Then*

$$Na_N^{-2}[d_{N,h}, h = 0, 1, \ldots, H]$$

$$\xrightarrow{d} \left(S - \frac{\alpha}{\alpha-2}\right) \left[\sum_{j=0}^{\infty} \psi(j)\psi(j+h), h = 0, 1, \ldots, H\right].$$

**Proof.** Observe that $D_N = N^{-1} \sum_{t=1}^{N} \sum_{j=0}^{\infty} c_j Z_{t-j}^2 - \sigma^2 \sum_{j=0}^{\infty} c_j$. Since the $c_j$ defined by (5.9) are absolutely summable, using Theorem 4.1 of Davis and Resnick (1985) and following the proof of Theorem 2.2 of Davis and Resnick (1986), we conclude that $Na_N^{-2}D_N \xrightarrow{d} (S - \frac{\alpha}{\alpha-2}) \sum_{j=0}^{\infty} c_j$. □

**Lemma 5.2.** *Suppose* $2 < \alpha < 4$ *and that (2.2) and (2.4) hold. If* $d < 1/\alpha$ *(and* $a_N^{-4}N \ln N \to 0$ *if* $d = 1/4$*), then* $Na_N^{-2}r_{N,h} \xrightarrow{P} 0$.

**Proof.** Set $\xi_t(h) = \sum_{i,j \geq 0, i \neq j+h} \psi(i)\psi(j)Z_{t-j}Z_{t+h-i}$ so that $r_{N,h} = N^{-1} \sum_{t=1}^{N} \xi_t(h)$. If $\mathbf{E}Z = 0$ and $\mathbf{E}Z^2 = \sigma^2 < \infty$, then for $0 < d < 1/2$,

$$\mathbf{E}[\xi_n \xi_0] = \sigma^4 \sum_{n \leq i \neq j < \infty} [\psi(i)\psi(j)\psi(i-n)\psi(j-n) + \psi(i)\psi(j)\psi(j-n)\psi(i-n)],$$

which implies $\mathbf{E}[\xi_n\xi_0] \sim K_d\sigma^4 n^{4d-2}$. Therefore,

$$\text{Var}\left[\sum_{t=1}^{N}\xi_t\right] = \begin{cases} O(N), & \text{if } 0 < d < 1/4, \\ O(N\ln N), & \text{if } d = 1/4, \\ O(N^{4d}), & \text{if } 1/4 < d < 1/2. \end{cases} \tag{5.10}$$

The claim follows because $a_N = N^{1/\alpha}L(N)$, where $L(\cdot)$ is slowly varying. □

The following lemma establishes asymptotics for the diagonal terms in the case of finite fourth moment.

Denote by $\mathbf{N}(0, \boldsymbol{\Sigma})$ an $(H+1)$-variate normal vector with the covariance matrix $\boldsymbol{\Sigma} = \{\Sigma(h, h'), h, h' = 0, 1, \ldots, H\}$ given by

$$\Sigma(h, h') = \mathbf{E}\varepsilon_0^2 \sum_{i=0}^{\infty}\sum_{k=0}^{\infty} c_i(h)c_k(h'), \tag{5.11}$$

in which the $c_j(h)$ are defined by (5.9).

**Lemma 5.3.** *Suppose that* $\mathbf{E}Z^4 < \infty$ *and* (2.4) *holds. Then*

$$N^{1/2}[d_{N,h}, h = 0, 1, \ldots, H] \xrightarrow{d} \mathbf{N}(0, \boldsymbol{\Sigma}).$$

**Proof.** Use Theorem 7.7.8 of Anderson (1971). □

In the remainder of this section, we study the off-diagonal terms in the case $d > 1/4$.

Assuming only $\mathbf{E}Z^2 < \infty$, using Lemma 4.1, we can define a sequence of bounded, independent, identically distributed random variables $Z_i(T)$ such that $\mathbf{E}Z_i(T) = 0$, $\mathbf{E}[Z_i(T) - Z_i] \to 0$ and $Z_i(T)$ is independent of $\{Z_j(T), j \neq i\}$. Further, define

$$R_N(T) = \sum_{h=0}^{H} u_h r_{N,h}(T),$$

where

$$r_{N,h}(T) = \frac{1}{N}\sum_{t=1}^{N}\sum_{i \neq j+h}\psi(j)\psi(i)Z_{t-j}(T)Z_{t+h-i}(T)$$

$$= \frac{1}{N}\sum_{k \neq k'}Z_k(T)Z_{k'}(T)\sum_{t=1}^{N}\psi(t-k)\psi(t+h-k').$$

Denote by $\{W(t)\}$ the standard Wiener process on the real line.

**Lemma 5.4.** *Suppose* $\mathbf{E}Z^2 < \infty$ *and* (2.4) *holds. If* $d > 1/4$*, then*

$$N^{1-2d}[r_{N,h}(T), h = 0, 1, \ldots, H]$$

$$\xrightarrow{d} \sigma^2(T) \left[ \int_{-\infty}^{\infty} \int_{-\infty}^{\infty} f(x_1, x_2) W(\mathrm{d}x_1) W(\mathrm{d}x_2), h = 0, 1, \ldots, H \right],$$

*where* $\sigma^2(T) = \mathrm{Var}[Z_0(T)]$ *and*

$$f(x_1, x_2) = C_d^2 \int_0^1 (v - x_1)_+^{d-1} (v - x_2)_+^{d-1} \, \mathrm{d}v. \qquad (5.12)$$

**Proof.** Introducing

$$\tilde{C}_N(k, k') = \sum_{h=0}^{H} u_h \sum_{t=1}^{N} N^{1-2d} \psi(t - k) \psi(t + h - k'),$$

we obtain

$$R_N(T) = N^{2d-2} \sum_{k \neq k'} \tilde{C}_N(k, k') Z_k(T) Z_{k'}(T) = N^{2d-1} \sigma^2(T) \sum_{k \neq k'} \tilde{C}_N(k, k') Z_k^*(T) Z_{k'}^*(T),$$

where $Z_k^*(T) = N^{-1/2} Z_k(T)/\sigma(T)$.

Clearly, $\mathbf{E}Z_k^*(T) = 0, \mathrm{Var}[Z_k^*(T)] = 1/N$, and by the functional central limit theorem, $\sum_{aN \leq k \leq bN} Z_k^*(T) \xrightarrow{d} W(b) - W(a)$. By Lemma 4.3,

$$\int_{-\infty}^{\infty} \int_{-\infty}^{\infty} [\tilde{f}_N(x_1, x_2) - \tilde{f}(x_1, x_2)]^2 \to 0,$$

where

$$\tilde{f}_N(x_1, x_2) = \tilde{C}_N(k, k') \mathbf{I}_{[k/N, (k+1)/N)}(x_1) \mathbf{I}_{[k'/N, (k'+1)/N)}(x_2)$$

and $\tilde{f}(x_1, x_2) = f(x_1, x_2) \sum_{h=0}^{H} u_h$.

Thus, all conditions of Lemma 4.2 are satisfied and the claim follows. $\qquad \square$

The next result shows that Lemma 5.4 remains valid without assuming bounded errors.

**Lemma 5.5.** *Suppose* $\mathbf{E}Z^2 < \infty$ *and* (2.4) *holds. If* $d > 1/4$ *(and* $d < 1/2$*), then*

$$N^{1-2d}[r_{N,h}, h = 0, 1, \ldots, H] \xrightarrow{d} \sigma^2 \left[ \int_{-\infty}^{\infty} \int_{-\infty}^{\infty} f(x_1, x_2) W(\mathrm{d}x_1) W(\mathrm{d}x_2), h = 0, 1, \ldots, H \right],$$

*where* $f(x_1, x_2)$ *is defined in* (5.12)*.*

**Proof.** Since $\sigma^2(T) \to \sigma^2$, it is enough to show that for any $\epsilon > 0$ and each $h = 0, 1, \ldots, H$,

$$\lim_{T \to \infty} \limsup_{N \to \infty} P\{N^{1-2d}|r_{N,h} - r_{N,h}(T)| > \epsilon\} = 0. \tag{5.13}$$

Let

$$\xi_{t,h} = \sum_{k \neq k'} \psi(t-k)\psi(t+h-k')Z_k Z_{k'}$$

and

$$\xi_{t,h}(T) = \sum_{k \neq k'} \psi(t-k)\psi(t+h-k')Z_k(T)Z_{k'}(T).$$

Relation (5.13) will follow once we have shown that

$$\mathbf{E}\left[\sum_{t=1}^{N}(\xi_{t,h} - \xi_{t,h}(T))\right]^2 = \mathbf{E}[Z - Z(T)]^2 O(N^{4d}). \tag{5.14}$$

Observe that

$$\left[\sum_{t=1}^{N}(\xi_{t,h} - \xi_{t,h}(T))\right]^2 \leq 2[S_1(N,T) + S_2(N,T)],$$

where

$$S_1(N,T) = \left[\sum_{t=1}^{N}\sum_{k \neq k'} \psi(t-k)\psi(t+h-k')Z_k(Z_{k'} - Z_{k'}(T))\right]^2,$$

$$S_2(N,T) = \left[\sum_{t=1}^{N}\sum_{k \neq k'} \psi(t-k)\psi(t+h-k')Z_{k'}(T)(Z_k - Z_k(T))\right]^2.$$

To compute the expected value

$$\mathbf{E}S_1(N,T) = \sum_{t,s=1}^{N}\sum_{k \neq k'}\sum_{i \neq i'} \psi(t-k)\psi(t+h-k')\psi(s-i)\psi(s+h-i')$$

$$\times \mathbf{E}[Z_k(Z_{k'} - Z_{k'}(T))Z_i(Z_{i'} - Z_{i'}(T))],$$

note that the expected value on the right-hand side vanishes, except in two cases: (a) $k = i$ and $k' = i'$; (b) $k = i'$ and $k' = i$.

In case (a),

$$\mathbf{E}[Z_k(Z_{k'} - Z_{k'}(T))Z_i(Z_{i'} - Z_{i'}(T))] = \sigma^2 \mathbf{E}[Z - Z(T)]^2.$$

In case (b), by Schwarz's inequality,

$$\mathbf{E}[Z_k(Z_{k'} - Z_{k'}(T))Z_i(Z_{i'} - Z_{i'}(T))] \leq \sigma^2 \mathbf{E}[Z - Z(T)]^2.$$

It follows that

$$\mathbf{E}S_1(N, T) \leq \sigma^2 \mathbf{E}[Z - Z(T)]^2[E_{11}(N, h) + E_{12}(N, h)],$$

where

$$E_{11}(N, h) = \sum_{s,t=1}^{N} \sum_{k \neq k'} \psi(t-k)\psi(t+h-k')\psi(s-k)\psi(s+h-k'),$$

$$E_{12}(N, h) = \sum_{s,t=1}^{N} \sum_{k \neq k'} \psi(t-k)\psi(t+h-k')\psi(s-k')\psi(s+h-k).$$

We must thus verify that $E_{11}(N, h)$ and $E_{12}(N, h)$ are $O(N^{4d})$. We will show the verification for $E_{12}(N, h)$ because the argument for $E_{11}(N, h)$ is the same (in fact, it is slightly shorter, as the $h$ cancels in the corresponding sums). Setting $i = t - k$ and $i' = t + h - k'$, we obtain

$$|E_{12}(N, h)| \leq \sum_{s,t=1}^{N} \sum_{i} |\psi(i)\psi(s-t+h+i)| \sum_{i'} |\psi(i')\psi(s-t-h+i')|.$$

We apply to the right-hand side the method of summing over the diagonals and consider the cases $s - t = 0, n = s - t > 0$ and $-n = s - t < 0$. Consequently,

$$|E_{12}(N, h)| \leq \bar{E}_{120} + \bar{E}_{12+} + \bar{E}_{12-},$$

where

$$\bar{E}_{120} = N \sum_{i} |\psi(i)\psi(h+i)| \sum_{i'} |\psi(i')\psi(-h+i')| \leq N\left[\sum_{i} \psi^2(i)\right]^2 = O(N),$$

$$\bar{E}_{12+} = \sum_{n=1}^{N-1} (N-n) \sum_{i} |\psi(i)\psi(n+h+i)| \sum_{i'} |\psi(i')\psi(n-h+i')|,$$

$$\bar{E}_{12-} = \sum_{n=1}^{N-1} (N-n) \sum_{i} |\psi(i)\psi(-n+h+i)| \sum_{i'} |\psi(i')\psi(-n-h+i')|.$$

A change of variables shows that $\bar{E}_{12-} = \bar{E}_{12+}$, so it remains to verify that $\bar{E}_{12+} = O(N^{4d})$. Since, by (2.4), $|\psi(j)|/j^{d-1}$ is bounded, we have, for any $n \geq 0$,

$$\sum_{i} |\psi(i)\psi(n+i)| = O\left(\int_{1}^{\infty} i^{d-1}(n+i)^{d-1}\,\mathrm{d}i\right) = O(n^{2d-1}).$$

It follows that

$$\bar{E}_{12+} = O(N)\left\{ \sum_{n=h+1}^{N-1} (n+h)^{2d-1}(n-h)^{2d-1} + \sum_{n=1}^{h-1} (n+h)^{2d-1}(h-n)^{2d-1} \right\}.$$

The second term in the braces is $O(h^{2d}) = O(1)$. The first term is

$$\sum_{j=1}^{N+h-1} (j+2h)^{2d-1} j^{2d-1} \leq \sum_{j=1}^{N+h-1} j^{4d-2} = O(N^{4d-1}).$$

This completes the verification that $\bar{E}_{12+} = O(N^{4d})$. We have thus established that $\mathbf{E}S_1(N,T) = \sigma^2 \mathbf{E}[Z - Z(T)]^2 O(N^{4d})$. Exactly the same argument applies to $S_2(N,T)$, so (5.14) and, consequently, (5.13) follow. □

# Acknowledgements

Lajos Horváth partially supported by NSF Grant DMS-06-04670. Piotr Kokoszka partially supported by NSF Grant DMS-04-13653.